\documentclass[review]{elsarticle}
\usepackage[latin1]{inputenc}
\usepackage[english]{babel}
\usepackage{amssymb,srcltx,amsthm,epsfig,bm}
\usepackage{amssymb,amsthm,epsfig}
\usepackage[centertags]{amsmath}
\usepackage{graphicx}
\usepackage{graphics}
\usepackage{natbib, rotating}
\usepackage{color}
\usepackage{booktabs, subfigure}
\usepackage{bbm}

\usepackage[paperwidth=210mm,paperheight=290mm,hmargin={25mm,25mm},vmargin={30mm,30mm}]{geometry}

\newcommand{\vt}{\mathbf{t}}

\newcommand{\vu}{\mathbf{u}}
\newcommand{\vA}{\mathbf{A}}
\newcommand{\vb}{\mathbf{b}}
\newcommand{\vZ}{\mathbf{Z}}
\newcommand{\vY}{\mathbf{Y}}
\newcommand{\vP}{\mathbf{P}}
\newcommand{\vI}{\mathbf{I}}

\newtheorem{thm}{Theorem}[section]
\newtheorem{lem}[thm]{Lemma}
\newtheorem{proposition}[thm]{Proposition}
\newcommand{\distras}

\newcommand{\vepsilon}{\mbox{\boldmath $\epsilon$}}

\newcommand{\pkonv}{\stackrel{p}{\rightarrow}}

\newcommand{\bqa}{\begin{eqnarray*}}
\newcommand{\eqa}{\end{eqnarray*}}
\newcommand{\bqan}{\begin{eqnarray}}
\newcommand{\eqan}{\end{eqnarray}}
\newcommand{\bit}{\begin{itemize}}
\newcommand{\eit}{\end{itemize}}
\newcommand{\ben}{\begin{enumerate}}
\newcommand{\een}{\end{enumerate}}
\newcommand{\beq}{\begin{equation}}
\newcommand{\eeq}{\end{equation}}
\newcommand{\bdes}{\begin{description}}
\newcommand{\edes}{\end{description}}
\date{}

\begin{document}

\begin{frontmatter}
\title{Consistent Variable Selection for Functional Regression Models}

\author{Julian A. A. Collazos}
\author{Ronaldo Dias}
\address{Department of Statistics - State University of Campinas (UNICAMP)}
\address{Rua Sergio Buarque de Holanda, 651, Distr. de Barao Geraldo, Campinas, Sao Paulo, Brazil}
\author{Adriano Z. Zambom\corref{mycorrespondingauthor}}
\address{Department of Statistics - Penn State University}
\address{323 Thomas Bldg., University Park, PA}

\cortext[mycorrespondingauthor]{Corresponding author}
\ead{adriano.zambom@gmail.com}

\begin{abstract}
The dual problem of testing the predictive significance of a particular covariate, and identification of the set of relevant covariates is common in applied research and methodological investigations. 
To study this problem in the context of functional linear regression models with predictor variables observed over a grid and a scalar response, 
we consider basis expansions of the functional covariates and apply the likelihood ratio test. Based on p-values from testing each predictor, we propose a new variable selection method, which is 
 consistent in selecting the relevant predictors from set of available predictors that is allowed to grow with the sample size $n$.
Numerical simulations suggest that the proposed variable selection procedure outperforms existing methods found in the literature. A real dataset from weather stations in Japan is analyzed.
\end{abstract}

\begin{keyword}
B-splines, hypotheses testing, False Discovery Rate, Functional Data, likelihood ratio test
\end{keyword}

\end{frontmatter}

\section{Introduction}

 \indent \indent In regression analysis, selecting the relevant set of predictors is a fundamental step for building a good predictive model. Including insignificant predictors results in over-complicated models with less predictive power and reduced ability to discern and interpret the influence of each variable. However, classical selection methods have to be adapted to the high-dimensional data sets which are becoming increasingly common in several areas of research.
 
 When the data is observed at several time (or space) points, simple linear regression models cannot be directly used. Functional regression models (FRM) express the discrete observations of the predictor as a smooth function, and inference can then be made about a response variable based on the functional data (Ramsay and Silverman, 2005). Such models have become increasingly useful due to their large number of applications, see Kokozsca and Horvath (2012) for some fundamental results and Ferraty and Vieu (2006) for a nonparametric approach. This high demand has recently leveraged important theoretical advances, see for example James (2002), Ferraty and Vieu (2009), James, Wang and Zhu (2009), Ferraty, Laksaci, Tadj and  Vieu (2010), and Aneiros and Vieu (2013), Goia and Vieu (2014), to cite a few.

However, only a few authors have considered variable selection in functional regression analysis.
Aneiros and Vieu (2014) show how to perform variable selection using the continuous structure of the functional predictors by studying which of the discrete observed points should be incorporated. Using a partial linear model for multi-functional data, Aneiros and Vieu (2015) propose a variable selection method based on the continuous specificity of the functional data. Cuevas (2014, Section 5) presents an interesting overview of recent methods for functional data analysis including functional regression. Most recent contributions in regression for these models can be found in Bongiorno et al. (2014).  
Another class of such methods uses regularization techniques, 
where the penalty simultaneously shrinks parameters and selects variables. Matsui and Konishi (2011) studied the group SCAD regularization for estimating and selecting functional regressors while Mingotti, Lillo and Romo (2013) and Hong and Lian (2011) generalized the Lasso for the case of scalar regressors and a functional response. 
Other recent contributions to the variable selection problem in functional models are Fan and Li (2004), Aneiros, Ferraty, and Vieu (2011), Gertheiss, Maity, and Staicu (2013) and Ma, Song and Wang (2013).

In this paper, we propose a different approach, exploiting the conceptual connection between model
testing and variable selection: dropping a covariate from the model is equivalent to not rejecting the null hypothesis that its corresponding parameter(s) is equal to zero. Abramovich, Benjamini, Donoho and Johnstone (2006) showed that the application of a false discovery rate (FDR) controlling procedure, such as Benjamini and Yekutieli (2001), on p-values resulting from testing each null hypothesis can be translated into minimizing a model selection criterion. The extension and adaptation of the theory of hypothesis testing to functional models have been studied by several authors in the literature (Cardot,  Goia, and Sarda, 2004, Yang and Nie, 2008, Swihart, Goldsmith and Crainiceanu, 2013, Kong, Staicu and Maity, 2013, McLean, Hooker and Ruppert, 2014, Pomann, Staicu and Ghosh, 2014). An interesting application can be found in Meinshausen, Meier and Buhlmann (2009), with results on the connection between p-values and variable selection in regression analysis.

The main objective of this paper is twofold: study the asymptotic properties of the hypothesis test based on residual sum of squares for the relevance of a predictor in a multivariate functional regression model; and propose a competitive variable selection procedure based on FDR (or Bonferroni) corrections applied on the p-values from the tests of each available functional predictor.
The proposed test statistic is a likelihood ratio type test, where restricted and full models are estimated through the B-Splines basis expansions of both coefficients and functional predictors. 
We examine the shift (non-centrality parameter) of the distribution of the test statistic under the alternative hypothesis, which provides insight into the power of the test and induce the demonstration of consistency of the variable selection procedure.

The remainder of this paper is as follows. In Section 2, we formally describe the regression model with functional covariates and scalar response via basis expansions. In Section 3, we present the testing procedure and the variable selection method.
In Section 4 we evaluate the finite sample performance of the proposed variable selection through simulation examples and a real application with weather data is considered in Section 5.

\vspace{-.3cm}
\section{The functional regression model: FRM}
\indent \indent Suppose that we have $n$ observations $\lbrace(y_{i},\boldsymbol{x}_{i}(\vt)): \vt \in \boldsymbol{\mathcal{T}}, i=1,...,n\rbrace$, where $y_{i}$ is a scalar response, $\boldsymbol{x}_{i}(\vt)=(x_{i1}(t_1),...,x_{iM}(t_M))$ are functional predictors and $\boldsymbol{\mathcal{T}} = \mathcal{T}_1 \times \ldots \times \mathcal{T}_M$. 
Each $\mathcal{T}_m, m = 1, \ldots, M$, is a compact set in $\mathbb{R}$ where the $m$-th predictor may be observed. The functional predictors $x_m, m = 1, \ldots, M$ are assumed to be in a fixed design so that in practice  $t_m \in \mathcal{T}_m$ is a grid representing time or space. 
Suppose that each of the $M$ functional predictors can be expressed as:
\bqan\label{func.pred}
x_{im}(t_m)=\sum_{j=1}^{p_{m}} \omega_{imj} \phi_{mj}(t_m) = \boldsymbol{W}_{im}^{T} \boldsymbol{\phi}_{m}(t_m), \hspace{.3cm} m = 1, \ldots, M, t_m \in \mathcal{T}_m,
\eqan
where $\boldsymbol{W}_{im}=(\omega_{im1},...,\omega_{imp_{m}})^{T}$ are the vectors of coefficients and $\boldsymbol{\phi}_{m}(t_m)=(\phi_{m1}(t_m),...,\phi_{mp_{m}}(t_m))^{T}$ are vectors of B-Splines basis functions. 
The basis functions and the $p_m$ coefficients in (\ref{func.pred}) are assumed to be determined prior to the regression modeling through smoothing methods.
In general this finite B-splines representation of a functional predictor is a good approximation of smooth functions, such as functions in the Sobolev Space (see Reif, 1997).

We consider the functional regression model (Ramsay and Silverman, 2005) given by
\begin{equation}
\label{FRM1}
y_{i}=\beta_{0}+\sum_{m=1}^{M}\int_{\mathcal{T}_m} x_{im}(t_m)\beta_{m}(t_m) dt_m + \varepsilon_{i}, 
\end{equation}
where $\beta_0$ is a constant, $\varepsilon_{i}, i = 1, \ldots, n$ are i.i.d. Gaussian noises with mean  $0$ and constant variance $\sigma^{2}$, and $\beta_{m}(t_m)$ are functional coefficients that we assume can be 
represented through the basis expansion
\bqan\label{func.par}
\beta_{m}(t_m)=\sum_{j=1}^{p_{m}} b_{mj}\phi_{mj}(t_m)=\boldsymbol{b}_{m}^{T}\boldsymbol{\phi}_{m}(t_m), \hspace{.3cm} m = 1, \ldots, M,  t_m \in \mathcal{T}_m,
\eqan
for the parameter vectors $\boldsymbol{b}_{m}=(b_{m1},...,b_{mp_{m}})^{T}$.
Thus the FRM in (\ref{FRM1}) can be re-expressed as a linear model in the following way
\bqa
y_{i} &=& \beta_{0}+\sum_{m=1}^{M} \int_{\mathcal{T}_m} \boldsymbol{W}_{im}^{T} \boldsymbol{\phi}_{m}(t_m) \boldsymbol{\phi}_{m}^{T}(t_m)\boldsymbol{b}_{m} dt_m + \varepsilon_{i} = \beta_{0}+\sum_{m=1}^{M} \boldsymbol{W}_{im}^{T} \int_{\mathcal{T}_m} \boldsymbol{\phi}_{m}(t_m) \boldsymbol{\phi}_{m}^{T}(t_m)dt_m \boldsymbol{b}_{m} + \varepsilon_{i}\\
&=& \beta_{0}+\sum_{m=1}^{M} \boldsymbol{W}_{im}^{T} \boldsymbol{J_{\phi_{m}}} \boldsymbol{b}_{m} + \varepsilon_{i}
= \boldsymbol{Z_{i}}^{T}\boldsymbol{b} + \varepsilon_{i}, 
\eqa
or in matrix form $\vY = \vZ \vb + \vepsilon$,
where $\boldsymbol{Z}_{i}=(1,\boldsymbol{W}_{i1}^{T} \boldsymbol{J_{\phi_{1}}},...,\boldsymbol{W}_{iM}^{T} \boldsymbol{J_{\phi_{M}}})^{T}$, $\boldsymbol{b}=(\beta_{0},\boldsymbol{b}_{1}^{T},...,\boldsymbol{b}_{M}^{T})^{T}$, $\vZ = (\vZ_1, \ldots, \vZ_n)^T$, $\boldsymbol{J_{\phi_{m}}}=\int_{\mathcal{T}_m} \boldsymbol{\phi}_{m}(t_m) \boldsymbol{\phi}_{m}^{T}(t_m)dt_m$ are $p_{m}\times p_{m}$ cross product matrices and $\vepsilon$ is the vector of error terms. 
Since we adopt B-splines basis expansions, the cross product matrix $\boldsymbol{J_{\phi_{m}}}$ can be easily computed using the procedure in Kayano and Konishi (2009).

\vspace{-.3cm}  
\section{Methodology} \label{sec:methodology}
\subsection{Testing procedure}
\indent \indent In this section we address the problem of testing the relevance of an individual functional predictor in the multivariate FRM. We consider testing the $r$-th ($r \in \{1,\ldots, M\}$) predictor through the following null hypothesis 
\begin{equation}
\label{hypoth}
H_{0}^r:\boldsymbol{b}_{r}=\boldsymbol{0} \:\:\:\: vs \:\:\:\: H_{a}^r:\boldsymbol{b}_{r}\neq \boldsymbol{0}.
\end{equation}
In linear models with normal errors, least squares estimates, which minimize the residual sum of squares, are equivalent to maximum likelihood estimates. For ease of notation, in this section, we omit from all statistics the index $r$ that identifies the predictor being tested.
Let $\zeta$ and $\Omega$ denote the spaces generated by the predictors under $H_{0}$ and $H_{a}$ respectively. Note that $\zeta\subset\Omega$ and hence rank$(\Omega) = 1 + \sum_{m=1}^M p_m := k$ and rank$(\zeta) = k-p_r = 1 + \sum_{m=1}^Mp_m - p_r :=k_{0}$. We assume throughout this paper that the matrix $\vZ$ has full rank, that is, Z has $k < n$ linearly independent columns (see also condition (C1) in Section \ref{sec.consist}). This assumption guarantees the existence and uniqueness of the least squares estimators. Let $RSS_{0}$ and $RSS$ denote the residual sum of squares under $H_{0}$ and $H_{a}$ respectively, that is,
\bqan \label{RSS}
RSS_{0}&=&\displaystyle\sum_{i=1}^{n}\left(y_{i}-\boldsymbol{Z}_{i}^{T}\hat{\boldsymbol{b}}^0\right)^2 \mbox{ and } RSS =\displaystyle\sum_{i=1}^{n}\left(y_{i}-\boldsymbol{Z}_{i}^{T}\hat{\boldsymbol{b}}\right)^2,
\eqan
where $\hat{\boldsymbol{b}}^0 = \hat{\boldsymbol{b}} - (\vZ^T\vZ)^{-1}\vA^T(\vA(\vZ^T\vZ)^{-1}\vA^T)^{-1}\vA\hat{\boldsymbol{b}}$ for a $p_r \times k$ matrix $\vA$ defining the null hypothesis, i.e., $\vA\boldsymbol{b} = \boldsymbol{0}$ implies $\boldsymbol{b}_{r}=\boldsymbol{0}$.

For insight into the distribution of the test statistic and the non-centrality parameter presented below, it is useful to express the sum of squares $RSS_{0}$ and $RSS$ as a quadratic form. 
We write $\hat{\vY}_{0}=\vZ\hat{\vb}^{0} = \vP_{0}\vY$ and $\hat{\vY} = \vZ\hat{\vb} = \vP\vY$, where $\vP_{0}$ and $\vP$ are the orthogonal projection matrices which project $\vY$ onto the spaces $\zeta$ and $\Omega$, respectively. We can then rewrite the residual sum of squares as
$RSS_{0} = \vY^{T}(\vI_{n}-\vP_{0})\vY$ and $RSS = \vY^{T}(\vI_{n}-\vP)\vY$, so that
$RSS_{0}-RSS=\vY^{T}(\vP-\vP_{0})\vY$.
Since 
\bqa
\frac{RSS_{0}}{\sigma^{2}}\:\:\overset{H_{0}}{\sim}\:\: \boldsymbol{\chi_{n-k_{0}}^{2}} \:\:\:\: \text{and} \:\:\:\:\: \frac{RSS}{\sigma^{2}}\:\:\overset{H_{0}}{\sim} \:\:\boldsymbol{\chi_{n-k}^{2}},
\eqa
in order to test $H_0$ in (\ref{hypoth}) we use the likelihood ratio statistic
\bqan\label{Tl}
T_{L} &=& -2Ln\left[ \frac{\tilde{L}_{0}}{\tilde{L}}\right] =  -2\left[-\frac{1}{2\tilde{\sigma}^{2}}RSS_{0}+\frac{1}{2\tilde{\sigma}^{2}}RSS\right] 
 = \dfrac{RSS_{0}-RSS}{\tilde{\sigma}^{2}} \overset{H_{0}}{\underset{n\rightarrow \infty}{\rightarrow}} \:\: \boldsymbol{\chi_{k-k_{0}}^{2}}
\eqan
in distribution, with $\tilde{\sigma}^2 = RSS/n \pkonv \sigma^2$ the maximum likelihood ratio statistic. 
From the Normality assumption of the residuals and the fact that
\bqa
 \frac{1}{\sigma^2}E\left[RSS_0 - RSS\right] &=& \frac{1}{\sigma^{2}}\left[\sigma^{2}Tr(\vP-\vP_{0})+(\vZ\vb)^{T}(\vP-\vP_{0})\vZ\vb\right]
= (k - k_{0})+\delta = p_r + \delta,
\eqa
where 
\bqan\label{delta}
\delta = \vb^{T}\vZ^{T}(\vP - \vP_{0})\vZ\vb/\sigma^{2},
\eqan
the following proposition can be established.
\begin{proposition} (Theorem 5.3c in Rencher and Schaalje, 2008) Let $RSS$ and $RSS_0$ be defined as in (\ref{RSS}). Then, under the alternative hypothesis in (\ref{hypoth})
\bqa
\frac{RSS_{0}}{\sigma^{2}}\:\:\overset{H_{a}}{\sim}\:\: \boldsymbol{\chi_{n-k_{0}}^{2}}(\delta) \:\:\:\: \text{and} \:\:\:\:\: \frac{RSS}{\sigma^{2}}\:\:\overset{H_{a}}{\sim} \:\:\boldsymbol{\chi_{n-k}^{2}}, \:\:\:\: \text{ so that }\:\:\:\: \dfrac{RSS_{0}-RSS}{\sigma^{2}}\:\:\overset{H_{a}}{\underset{}{\sim}} \:\: \boldsymbol{\chi_{k - k_{0}}^{2}}(\delta).
\eqa
\end{proposition}

Lemma \ref{lemma 3.1} specifies the order of the non-centrality parameter of the distribution of $(RSS_0 - RSS)/\sigma^2$. Growing at the order of the sample size, multiplied by the significance size of the parameter being tested, the shift produced by the non-centrality parameter under $H_a$ provides evidence for rejecting the null hypothesis. Using this result, Theorem \ref{theorem 3.4} shows the consistency of the proposed variable selection procedure, which is described in Section \ref{sec.consist}. 

\begin{lem} \label{lemma 3.1} Let $T_{L}$ be the likelihood ratio test statistic defined in (\ref{Tl}) for testing $H_{0}$ in (\ref{hypoth}). For the alternative hypothesis,
the non-centrality parameter $\delta$ 
defined in (\ref{delta}) is of order $\delta \sim c (n - k_0)$, for a constant $c$.

\end{lem}

\vspace{-.3cm}
\subsection{Consistent test based variable selection} \label{sec.consist}
\indent \indent In this section we describe a test-based variable selection method which is shown to consistently identify the set of relevant predictors. A similar procedure was used by Bunea, Wegkamp and Auguste  (2006) in the linear model setting, and by Zambom and Akritas, (2014) for a nonparametric model. 

Let $I_{M}=\left\lbrace 1,...,M\right\rbrace$ denote the set of indices of the $M$ available functional predictors. Assume that the true underlying model is sparse in the sense that only a few predictors significantly relate to the response variable, while $M$ is allowed to grow with $n$ at a rate such that the following condition holds
\bqa
{\bf \mbox{Condition } (C1)}: \hspace{.3cm}   k = 1 +  \sum_{m=1}^M p_m \leq \sqrt{n}/log(n). 
\eqa
Let $I_{0}=\left\lbrace m_{1},...,m_{M_{0}}\right\rbrace$ denote the (unknown) subset of indices corresponding to the $M_{0}$ significant predictors. The objective of the proposed variable selection method is to identify the subset $I_{0}$, that is, to determine the set of functional variables with predictive significance.

Let $T_{L}^r$, $r=1,...,M$, denote the likelihood test statistic defined in (\ref{Tl}) for testing $H_{0}^{r}$ in (\ref{hypoth}) and 
\bqan\label{p_value}
\pi_{r}=1-\Psi(T_{L}^r)
\eqan
the corresponding p-value, where $\Psi(.)$ is the cumulative function of the $\boldsymbol{\chi_{p_r}^{2}}$ distribution. 
The Bonferroni method yields $\hat{I} = \{m: \pi_m \leq q/M\}$ as the estimate of $I_0$.
The false discovery rate (FDR) procedure (Benjamini and Yekutieli, 2001) computes
\begin{equation}
\label{FDR}
s=\max\left\lbrace j:\pi_{(j)}\leq\frac{j}{M}\frac{q}{\sum_{l=1}^{M}l^{-1}} \right\rbrace,
\end{equation}
where $\pi_{(1)}\leq ... \leq \pi_{(M)}$ denote the ordered p-values and $q$ is the choice of level, and rejects $H_{0}^{(j)}$, $j=1,...,s$. If no such $s$ exists, no hypothesis is rejected. The proposed variable selection method selects the predictors with indices corresponding to the $s$ rejected null hypotheses. Hence, $I_{0}$ is estimated by the set $\widehat{I}$ of indices corresponding to the first $s$ ordered p-values.

Let us now prove the consistency of the proposed variable selection method. Let $R$ denote the total number of rejected hypothesis, so we have that $R = s \mathbbm{1}(s \mbox{ in (\ref{FDR}) exists})$, where $\mathbbm{1}(.)$ is the indicator function.
Now, let $V$ be the number of falsely rejected hypotheses, and set $Q = (V/R)\mathbbm{1}(R > 0)$
for the proportion of falsely rejected hypotheses. By definition, the FDR is $E(Q)$, and $E(Q)\leq q(M-M_{0})/M\leq q$, (Benjamini and Yekutieli, 2001). We consider consistent a procedure, and the estimated set $\hat{I}$, if $P(\hat{I} = I_0) \rightarrow 1$ as $n \rightarrow \infty$.
Theorem \ref{theorem 3.4}, in connection with Lemmas \ref{lemma 3.1} - \ref{lemma 3.3}, show the consistency of $\hat{I}$.

\begin{lem} \label{lemma 3.2} Let $T_{L}^r$ and $\pi_{r}=1-\Psi(T_{L}^r)$  be the test statistic and the p-value defined as in (\ref{Tl}) and (\ref{p_value}) for testing $H_{0}^{r}$. Assume condition (C1) holds and define $A_n = \{|\tilde{\sigma} - \sigma| \leq \sqrt{log(n)/n}\}$. 

\begin{itemize}
\item[\textit{(a)}] \textit{For $r\notin I_{0}$ and any $0< \gamma < 1$, we have $P\left(\{\pi_{r}\leq\gamma\}\cap A_n\right)=\gamma+O(\sqrt{log(n)/n})$.}
\item[\textit{(b)}] \textit{For $r\in I_{0}$ and $0 < \gamma < 1$, as $n \rightarrow \infty$, if 
$\gamma \geq 1/n$, we have\\ $P\left(\{\pi_{r}>\gamma\}\cap A_n\right)=o(\gamma) + O(\sqrt{log(n)/n})$.}
\end{itemize}
\end{lem}

\begin{lem} \label{lemma 3.3} Let $\Gamma_{n}$ be the event where the smallest $M_{0}$ p-values defined in (\ref{p_value}) are the p-values corresponding to the $M_{0}$ significant functional predictors, with $I_{0}=\left\lbrace m_{1},...,m_{M_{0}}\right\rbrace$, that is
$$\Gamma_{n} = \left[ \left\lbrace \pi_{(1)},...,\pi_{(M_{0})}\right\rbrace =\left\lbrace \pi_{m_{1}},...,\pi_{m_{M_{0}}} \right\rbrace \right].$$
Then, if 
condition (C1) holds, $\lim\limits_{n\rightarrow \infty} P\left( \Gamma_{n}\right)=1$.
\end{lem}

\begin{thm} \label{theorem 3.4} Let $\delta$ be the non-centrality parameter defined in (\ref{delta}), and $q$ the chosen bound of FDR in (\ref{FDR}) or in Bonferroni corrections. Assume that 
condition (C1) holds and $q\rightarrow 0$ as $n\rightarrow\infty$, in such a way that $q \geq M\left(\sum_{l=1}^{M}l^{-1}\right)/(M_{0}n)$ and $Mq/log(M) \rightarrow 0$.
\textit{Then,} $\lim\limits_{n\rightarrow\infty}P\left(\hat{I}=I_{0}\right)=1.$
\end{thm}

Note that the choice of $q \rightarrow 0$ is important for the consistency of the proposed method. For real datasets, a rule of thumb is to choose $q = O(1/M)$ if $M$ is large relatively to the sample size $n$, otherwise  choose $q = O(1/\sqrt{n})$. These choices guarantee the consistency of the variable selection while satisfying all assumptions and conditions. In the simulation study we explore different choices of this parameter.

\vspace{-.3cm}
\section{Numerical simulations}
\indent \indent Simulation studies were conducted to evaluate the finite sample performance of the proposed variable selection procedure. 
The Monte Carlo simulations in this section are based on 100 and 300 generated observations of six functional covariates and a scalar response $\{(x_{im}(t), y_{i}); t \in \tau_{m}, i=1,..., n, m=1,..., 6\}$, 
 extending the simulation set up in Matsui and Konishi (2011) by including three extra functional predictors. We compared the performance of the proposed variable selection procedure with that of group SCAD and group LASSO proposed by Matsui and Konishi (2011), and the Generalized Functional Linear Model (GFLM) method in Gertheiss, et al. (2013) with adaptive penalization. For comparison purposes, we used 6 basis functions for the estimation of the predictors and the functional parameters $\beta(.)$ in all methods. First, we generated $z_{im}$ corresponding to the predictor $X_m$ in an equally spaced grid of 50 points in $\mathcal{T}_m$ in the following way:
\bqa
z_{im}=u_{im}(t_{m})+\epsilon_{im}, \:\:\:\:\: \epsilon_{im} \sim N(0, (0.025r_{x_{im}})^{2}),
\eqa
where $r_{x_{im}}=\max_{i}(u_{im}(t_{m}))-\min(u_{im}(t_{m}))$ and 

$u_{i1}(t) = cos(2\pi(t-a_{1}))+a_{2}, \:\:\:\:\: \mathcal{T}_{1}=[0, 1], \:\:\: a_{1}\sim N(-4, 3^2), \:\:\: a_{2}\sim N(7, 1.5^2),$

$u_{i2}(t)=b_{1}sin(\pi t)+b_{2}, \:\:\:\:\: \mathcal{T}_{2}=[0, \pi/3], \:\:\: b_{1}\sim U(3, 7), \:\:\: b_{2}\sim N(0, 1),$

$u_{i3}(t)=c_{1}t^{3} + c_{2}t^{2} + c_{3}t, \:\:\:\:\: \mathcal{T}_{3}=[-1, 1], \:\:\: c_{1}\sim N(-3, 1.2^2), \:\:\: c_{2}\sim N(2, 0.5^2), c_{3}\sim N(-2, 1),$

$u_{i4}(t)=sin(2(t-d_{1}))+d_{2}t, \:\:\:\:\: \mathcal{T}_{4}=[0, \pi/3], \:\:\: d_{1}\sim N(-2, 1), \:\:\: d_{2}\sim N(3, 1.5^2),$

$u_{i5}(t)=e_{1}cos(2t)+e_{2}t, \:\:\:\:\: \mathcal{T}_{5}=[-2, 1], \:\:\: e_{1}\sim U(2, 7), \:\:\: e_{2}\sim N(2, 0.4^2),$

$u_{i6}(t)= f_{1}e^{-t/3} + f_{2}t + f_{3}, \:\:\:\:\: \mathcal{T}_{6}=[-1, 1], \:\:\: f_{1}\sim N(4, 2^2), \:\:\: f_{2}\sim N(-3, 0.5^2), f_{3}\sim N(1, 1).$\\
The scalar response $y_{i}$ was generated as $y_{i}=g(\vu_{i})+\varepsilon_{i}$, 
where $g(\vu_{i})=\displaystyle\sum_{m=1}^{6}\int_{\mathcal{T}_{m}}u_{im}(t)\beta_{m}(t)dt$, $\varepsilon_{i}\sim N(0, (0.05R_{y_{i}})^{2})$ and $R_{y_{i}}=max(g(\vu_{i}))-min(g(\vu_{i}))$. For a constant $c = 0, 0.4$ and $0.8$, the coefficient functions $\beta_{m}(t)$ are given by
\bqa
\beta_{1}(t)=sin(t), \:\:\:\:\: \beta_{2}(t)=sin(2t), \:\:\:\:\: \beta_{3}(t)=-ct^2, \:\:\:\:\: \beta_{4}(t)=sin(2t), \:\:\: \beta_{5}(t)=csin(\pi t), \:\:\: \beta_{6}(t)=0.
\eqa
Note that if $c = 0$ the true model specifies that only $u_1, u_2$ and $u_4$ significantly relate to the response, corresponding to the predictors $X_1, X_2$ and $X_4$.

As the first step of our analysis, the random data $z_{im}$ was converted into the functional data $x_{im}$ using B-splines basis smoothing.
For these data, we assumed the functional regression model
\bqa
y_i = \sum_{m=1}^6 \int_{\mathcal{T}_{m}}x_{im}(t)\beta_{m}(t)dt + \varepsilon_i,
\eqa
and applied the proposed variable selection method described in Section \ref{sec:methodology}.
With 100 Monte Carlo simulations, we computed the number of correctly selected models and the averages of the mean square errors (AMSE) for the proposed method with FDR and Bonferroni corrections, as well as for group LASSO, group SCAD and GFLM. 
\begin{table}[!htb]
\footnotesize
\label{tabela}
\caption{Number of correctly selected models and AMSE} \label{tabela1}
\centering
\begin{tabular}{cllccccccccccc} 
\toprule
& & & \multicolumn{3}{c}{$T_{L}^{BC}$}& \multicolumn{3}{c}{$T_{L}^{FDR}$}  & \multicolumn{2}{c}{SCAD} &  \multicolumn{2}{c}{LASSO}& GFLM\\
 \cmidrule(lr){4-6}
 \cmidrule(lr){7-9}
 \cmidrule(lr){10-11}
 \cmidrule(lr){12-13}
$c$ & $n$ & &$.01$ & $.05$ & $.1$ & $.01$ & $.05$ & $.1$ & GCV & BIC & GCV & BIC \\
\midrule
0     & 100 &correct& 88 & 79 & 65 & 87 & 74 & 58 & 82 & 82 & 80 & 83 & 77\\
     &&AMSE&(2.07)&(2.04)&(2.01)&(2.06)&(2.05)&(1.97)&(1.45)&(1.45)&(1.19)&(1.30) & (8.94)\\ 
      & 300 &correct& 96 & 92 & 88 & 95 & 89 & 83 & 85 & 85 & 84 & 86 & 83\\
     &&AMSE&(1.93)&(1.98)&(1.89)&(1.92)&(1.97)&(1.91)&(1.31)&(1.31)&(1.04)&(1.16) & (8.51)\\
\hline
.4         & 100 &correct& 79 & 79 & 78 & 82 & 80 & 73 & 79 & 79 & 65 & 65 & 76\\
     &&AMSE&(2.61)&(2.98)&(2.77)&(2.88)&(3.01)&(2.82)&(5.60)&(5.60)&(5.67)&(5.70) & (11.37)\\
      & 300 &correct& 96 & 94 & 90 & 95 & 92 & 88 & 83 & 83 & 71 & 80 & 84\\
     &&AMSE&(2.57)&(2.90)&(2.74)&(2.87)&(2.91)&(2.79)&(5.58)&(5.58)&(5.64)&(5.59) & (10.78)\\
\hline
.8         & 100 &correct& 83 & 81 & 80 & 83 & 81 & 79 & 83 & 83 & 72 & 74 & 83\\
     &&AMSE&(7.15)&(7.96)&(7.92)&(7.42)&(7.87)&(7.78)&(7.41)&(7.41)&(7.14)&(7.87) & (13.49)\\
      & 300 &correct& 98 & 96 & 93 & 99 & 95 & 92 & 93 & 93 & 80 & 82 & 94\\
     &&AMSE&(7.08)&(7.10)&(7.01)&(7.09)&(7.11)&(7.14)&(7.27)&(7.27)&(7.17)&(7.32) & (12.05)\\      
\bottomrule
\end{tabular}
\end{table}
\normalsize
The results in Table \ref{tabela1} suggest that when the sample size is relatively small (n = 100), all four methods seem to select the correct model about the same number of times, however as the sample size increases, the proposed variable selection procedure outperforms group SCAD, group LASSO and the GFLM. We note that restrictive choices of level for the tests tend to yield better results of the proposed method, where for example we observe that the choice of $q = 0.01$ delivers the highest number of correctly model selections.
For $c = 0$ or $c = 0.8$, group SCAD and group LASSO have AMSE similar to that of the proposed procedure. However for predictors included in the model with low significance ($c = 0.4$), the AMSE of group SCAD and group LASSO are about double the AMSE achieved by our procedure, while the GFLM delivers the highest AMSE in all models.

\vspace{-.3cm}
\section{Real Data Example: Weather Data}
\indent \indent In this application, we consider weather data observed monthly at 79 weather stations in Japan. The data set was obtained from http://www.data.jma.go.jp/obd/stats/data/en/, and includes
monthly and annual total observations averaged from 1971 to 2000: monthly
observed average temperatures (TEMP), average atmospheric pressure (PRESS), time of daylight (LIGHT), average humidity
(HUMID), maximum temperature (MAX.TEMP), minimum temperature (MIN.TEMP) and annual total precipitation. The dataset used in this analysis does not correspond to the one used in Matsui and Konishi (2011), rather we selected the 79 most reliable stations according to the aforementioned website.

The functional predictors, observed at a grid of 1 to 12 points, were fitted using 6 B-splines basis functions. Figure \ref{weather1} shows examples of the fitted functional predictors.
The goal of this application is to select the functional covariates that significantly relate to annual total precipitation. 
We applied the proposed variable selection method and compared the results with those of the group SCAD, group LASSO and GFLM selection procedures, using the same number of basis functions.

\begin{center} 
[Figure 1 about here]
\end{center} 
\begin{figure}[!htb]
\vspace{-.8cm}
\caption{Examples of smoothed functional covariates from weather data}
\label{weather1}
\end{figure}

The selected functional predictors for each method are shown in Table \ref{tab_weather1}. Humidity and maximum temperature are selected by all methods except GFLM, however, differently from group SCAD and group LASSO, the proposed procedure and GFLM selected PRESS and did not select LIGHT. Atmospheric pressure is well known among meteorologists to be related to precipitation.
Low and high air pressure systems are usually caused by unequal heating across the surface of the planet. A low pressure system is an area where the atmospheric pressure is lower than that of the area around it. The production of clouds and consequent precipitation are hence related to the wind, warm air and atmospheric lifting caused by low pressure systems.
\begin{table}[!htb]
\small
\caption{Selected predictors for the weather dataset example} \label{tab_weather1}
\centering
\begin{tabular}{c  c} 
\toprule
Method & Selected \\
\hline
$T_L$ & PRESS, HUM, MAX.T\\
SCAD & LIGHT, HUM, MAX.T\\
LASSO & TEMP, LIGHT, HUM, MAX.T\\
GFLM & TEMP, PRESS, LIGHT\\
\bottomrule
\end{tabular}
\end{table}

In a simulation of 100 bootstrap samples from the weather data, we performed variable selection using the proposed method, group SCAD and group LASSO and GFLM. Table \ref{tab_weather2} shows the number of times each predictor was selected. While LIGHT was the third most selected predictor by group SCAD and group LASSO (about 70\% of the time) and the most selected by GFLM, it was only the fourth most selected predictor when using the proposed procedure. On the other hand, pressure was selected most frequently by the proposed method, followed by humidity and maximum temperature. 
Our results meet the expectations of most specialized meteorology literature, which finds significant relation between pressure, humidity and maximum temperature with annual precipitation.

\begin{table}[!htb]
\small
\caption{Ratio of selection on 100 bootstrap samples of weather data} \label{tab_weather2}
\centering
\begin{tabular}{ccccccc} 
\toprule
Method & TEMP &PRESS & LIGHT & HUM & MAX.T & MIN.T\\
\midrule
$T_{L}{(BC)}$ &0.38& 0.90& 0.56 & 0.89& 0.87& 0.41\\
$T_{L} {(FDR)}$ &0.40& 0.90& 0.58 & 0.87& 0.86& 0.45\\
\hline   
SCAD (GCV) &0.37&  0.23& 0.65&  0.81& 0.81&  0.24\\

SCAD (BIC) &0.37&  0.21& 0.75&  0.81& 0.83&  0.23\\
\hline
LASSO (GCV)&0.45& 0.35& 0.62&  0.78& 0.80& 0.25\\

LASSO (BIC)&0.45& 0.34& 0.75&  0.81& 0.81& 0.23\\
\hline
GLM & 0.73 & 0.67 & 0.79 & 0.47 & 0.47 & 0.21\\
\bottomrule
\end{tabular}
\end{table}

\vspace{.31cm}
{\bf Acknowledgements}
This paper was partially supported by CNPq (grant 302956/2013-1), Fapesp (grant 2013/07375-0 and 2013/00506-1) and CAPES. We would also like to thank Michael G. Akritas and Nancy Lopes Garcia for their fruitful insights.

{\large Appendix}

\vspace{.2cm}
{\it Proof of Lemma \ref{lemma 3.1}}\\
Since $(\vP-\vP_{0})$ is idempotent, it is easy to show that the non-centrality parameter $\delta$ is equal to
\bqa
\delta = \vb^{T}\vZ^{T}(\vP-\vP_{0})\vZ\vb/\sigma^{2} = ||\vZ\vb - \vP_0\vZ\vb||^2/\sigma^2.
\eqa
Note that $E(\vY|\vZ) = \vZ\vb$ is the vector of expected values conditional on $\vZ$, which belongs to the subspace $\Omega$, and $\vP_0\vZ\vb$ is its projection onto the restricted subspace $\zeta$.
Without loss of generality write $\vZ\vb = (\vZ^0, \vZ^1)(\vb_{-r},\vb_r)$, where $\vZ^1$ is the sub-matrix of $\vZ$ with columns corresponding to the parameters $\vb_r$, and $\vZ^0$ the remaining columns (similarly for $\vb_{-r}$).
Let $\tilde{\vY} = \vZ\vb$ so that $(\vP - \vP_0)\vZ \vb = \tilde{\vY} - \vP_0\tilde{\vY} = (\vI - \vP_0)\tilde{\vY}$. The quantity $(\vI - \vP_0)\tilde{\vY}$ is the residuals from the projection of $\tilde{\vY}$ onto the subspace $\zeta$. 
This can be viewed as a linear model $\tilde{\vY} = E(\tilde{\vY}|\vZ^0) + \tilde{\varepsilon}$, so that the mean squared error $||(\vI - \vP_0)\tilde{\vY}||^2/(n-k_0) = \tilde{\vY}^T(\vI - \vP_0)\tilde{\vY}/(n-k_0) = \delta\sigma^2/(n-k_0)$ will converge to the constant. This implies that $\delta \sim c(n - k_0)$.




\vspace{.2cm}
{\it Proof of Lemma \ref{lemma 3.2}}\\
{\it Part (a)} Let $\Psi_{p_r}(.)$ be the cumulative distribution function (c.d.f.) of the central $\chi^{2}_{p_r}$ distribution and $\Psi_{p_r}^{-1}(.)$ its inverse. Also, denote the residual sum of squares under hypothesis $H_0^r$ in (\ref{hypoth}) by $RSS_0^r$. Using the fact that $\lim_{n\rightarrow\infty}P(A_n) = 1$ (Lemma A.1 in Bunea et al., 2006), we obtain $\lim_{n\rightarrow\infty}P(|\tilde{\sigma}^2 - \sigma^2| \geq \sigma\alpha) = 0$ for $\alpha = \sqrt{log(n)/n}$. For all $r\notin I_{0}$, $\vb_{r}=0$, and for any $0 < \gamma < 1$ we find that 
\bqa
P\left(\left\{\pi_{r}\leq\gamma\right\}\cap A_n\right)=P\left(\left\{1-\Psi_{p_r}(T_{L}^r)\leq\gamma\right\}\cap A_n\right)
=P\left(\left\{T_{L}^r\geq \Psi_{p_r}^{-1}(1-\gamma)\right\}\cap A_n\right)\hspace{3.3cm}\\
=P\left(\left\{\frac{RSS_{0}^r-RSS}{\tilde{\sigma}^{2}}\geq \Psi_{p_r}^{-1}(1-\gamma)\right\}\cap A_n\right)
\leq P\left(\frac{RSS_{0}^r-RSS}{\sigma^{2}}\geq \left(1 - \frac{\alpha}{\sigma}\right)\Psi_{p_r}^{-1}(1-\gamma)\right)
=\gamma+O(\alpha).
\eqa

{\it Part (b)} Let $\alpha = \sqrt{log(n)/n}$. For all $0<\gamma < 1$,
\bqa
&&P\left(\left\{\pi_{r}>\gamma\right\}\cap A_n\right)=P\left(\left\{1-\Psi_{p_r}(T_{L}^r)>\gamma\right\}\cap A_n\right) =P\left(\left\{T_{L}^r< \Psi_{p_r}^{-1}(1-\gamma)\right\}\cap A_n\right)\\
&&=P\left(\left\{\frac{RSS_{0}^r-RSS}{\tilde{\sigma}^{2}} < \Psi_{p_r}^{-1}(1-\gamma)\right\}\cap A_n\right)\leq P\left(\frac{RSS_{0}^r-RSS}{\sigma^{2}} < \left(\frac{\alpha}{\sigma} + 1\right)\Psi_{p_r}^{-1}(1-\gamma)\right).
\eqa

Under the alternative $(RSS_{0}^r-RSS)/\sigma^{2}$ has a non-central chi-square distribution with $p_r$ degrees of freedom and non-centrality parameter $\delta$, whose c.d.f. we denote by $\Psi_{p_r, \delta}(.)$. 
Since $\delta \sim c(n - k_0)$ and $k \leq \sqrt{n}/log(n)$, we conservatively have $\delta \sim c(n - \sqrt{n}/log(n))$.
For $\gamma \geq 1/n$, as $n \rightarrow \infty$ and hence $\delta \rightarrow \infty$, we have that
\bqa
\Psi_{p_r, \delta}(\Psi_{p_r}^{-1}(1-\gamma)) &=& \sum_{j = 0}^\infty\frac{\delta^j}{2^j j!}e^{-\frac{\delta}{2}}\Psi_{p_r+2j}(\Psi_{p_r}^{-1}(1-\gamma))\\
 &=& \sum_{j = 0}^\infty\frac{\delta^j}{2^j j!}e^{-\frac{\delta}{2}}\left(1 - e^{-\Psi_{p_r}^{-1}(1-\gamma)/2}\sum_{\ell = 0}^{p_r/2+j-1}\frac{\left(\Psi_{p_r}^{-1}(1-\gamma)\right)^\ell}{2^\ell j!}\right)
= o(\gamma),
\eqa
since 
the poisson weights are dislocated to larger values of $j$ at a rate of $\exp(n - \sqrt{n}/log(n))$ while the values of $\Psi_{p_r+2j}(\Psi_{p_r}^{-1}(1-\gamma))$ are dislocated at a rate slower than $n$, for the choice of $\gamma$ (Note that even if $\gamma$ was chosen to decrease at a slower rate than $\exp(-n)n^k$, the percentile $\Psi_{p_r}^{-1}(1-\gamma)$ would increase slower than a linear rate in $n$, and $\Psi_{p_r, \delta}(\Psi_{p_r}^{-1}(1-\gamma))$ would be $o(1)$). 
Hence $P\left(\left\{\pi_{r}>\gamma\right\}\cap A_n\right) \leq \Psi_{p_r, \delta}(\left(\frac{\alpha}{\sigma} + 1\right)\Psi_{p_r}^{-1}(1-\gamma)) = o\left(\gamma\right) + O(\alpha)$.\hspace{1cm} $\square$

\vspace{.2cm}
{\it Proof of Lemma \ref{lemma 3.3}}\\
Since $\lim_{n\rightarrow\infty}P(A_n)  = \lim_{n\rightarrow \infty}P\left(|\tilde{\sigma} - \sigma| \leq \alpha\right)= 1$, where $\alpha = \sqrt{log(n)/n}$, it suffices to show that $\lim_{n\rightarrow\infty} P(\Gamma_{n}^{c} \cap A_n) = 0$.
From Lemma \ref{lemma 3.1}, $\delta$ is of order $\sim cn$, so that 
for $\gamma = \alpha$
$$
\begin{aligned}
P\left(\Gamma_{n}^{c} \cap A_n\right)&\leq \displaystyle\sum_{m\in I_{0}}\displaystyle\sum_{k\notin I_{0}}P\left(\{\pi_{k}<\pi_{m}\}\cap A_n\right)  \\
&\leq\displaystyle\sum_{m\in I_{0}}\displaystyle\sum_{k\notin I_{0}}\left[ P\left(\{\pi_{k}\leq\gamma\}\cap A_n\right)+P\left(\{\pi_{m}>\gamma\}\cap A_n\right) \right]\\
&\leq\displaystyle\sum_{m\in I_{0}}\displaystyle\sum_{k\notin I_{0}}\left[ \gamma+O(\alpha)+o(\gamma) \right] =M_{0}(M-M_{0})\left[ \gamma+O(\alpha)+o(\gamma) \right], 
\end{aligned}
$$
where the last inequality follows from Lemma \ref{lemma 3.2}.
Since $\gamma = \alpha$ we have $\lim\limits_{n\rightarrow \infty} P\left(\Gamma_{n}^c\cap A_n\right)=0$.\hspace{.6cm}$\square$

\vspace{.2cm}
{\it Proof of Theorem \ref{theorem 3.4}}\\
We follow the proof in Bunea et. al. (2006) to prove the theorem under FDR corrections. The case of Bonferroni corrections follows with similar steps.
If $\hat{I}$ is equal to $I_{0}$, we have $M_{0}$ rejections ($R=M_{0}$) with none of them being erroneous ($V=0$). Thus, the consistency of $\hat{I}$ is verified by showing that
\begin{equation}
P(\hat{I}=I_{0})=P(R=M_{0},V=0)\rightarrow 1,\:\: \text{as}\:\: n\rightarrow\infty.
\end{equation}
This follows by showing that both $P(R\neq M_{0})$ and $P(V\geq 1)$ are asymptotically negligible. We have that (Bunea et al. 2006, Lemma 2.1)
\begin{equation}
\label{buneal2.1}
P(V\geq 1)\leq P(R\neq M_{0})+\frac{M_{0}(M-M_{0})}{M}q.
\end{equation}

Hence, in order to show consistency of $\hat{I}$ we need only show that $P(R\neq M_{0})\rightarrow 0$. Let $q_{M}=q/ \sum_{l=1}^{M}l^{-1}$ and note that
$
\left\{R \neq M_0\right\} = \overset{M}{\underset{m = M_0+1}{\cup}}\left\{\pi(m) \leq q_M m/M\right\}\cup \left\{\pi(M_0) > q_M M_0/M\right\},
$
so that 
\bqan
\label{consist}
P( R\neq M_{0}) &\leq& P(A_n^c) + P(\Gamma^c\cap A_n) + P\left(\left\{\pi_{(M_{0})} > q_{M}\frac{M_{0}}{M}\right\}\cap\Gamma_n\cap A_n\right)\nonumber\\
&&+\sum_{m=M_{0}+1}^{M} P\left(\left\{\pi_{(m)}\leq q_{M}\frac{m}{M}\right\}\cap\Gamma_n\cap A_n\right),  
\eqan
where $A_n = \left\{|\tilde{\sigma} - \sigma| \leq \alpha\right\}$, with $\alpha = \sqrt{log(n)/n}$, and $\Gamma_n$ is the event defined in Lemma \ref{lemma 3.3}. The third term on the right hand side of (\ref{consist}) is equal to
\bqa
P\left(\left\{\pi_{(M_{0})} > q_{M}\frac{M_{0}}{M}\right\}\cap\Gamma_n\cap A_n\right)  &\leq& M_{0}\max_{m\in I_{0}} P\left(\left\{\pi_{m} > q_{M}\frac{M_{0}}{M}\right\}\cap A_n\right)\\
&& =O\left(M_0\left(o\left(\frac{q_M M_0}{M}\right) + \alpha \right)\right) =o(1), \:\: \text{as} \:\: n\rightarrow\infty,   
\eqa
by Lemma \ref{lemma 3.2} 
and the assumptions of the theorem. For the last term in (\ref{consist}) we have
\bqa
&&\sum_{m=M_{0}+1}^{M} P\left(\left\{\pi_{(m)}\leq q_{M}\frac{m}{M}\right\}\cap\Gamma_n\cap A_n\right) \leq \sum_{m=M_{0}+1}^{M}P\left(\left\{\pi_{(m)}\leq q_{M}\right\}\cap \Gamma_{n}\cap A_n\right)  \\
&&\hspace{1cm}\leq \sum_{m\notin I_{0}} P\left(\left\{\pi_{m}\leq q_{M}\right\}\cap A_n\right)
= O\left((M - M_0)\left(\frac{q}{log(M)} + \alpha\right)\right) = o(1), \:\: \text{as} \:\: n\rightarrow\infty, 
\eqa
by Lemma \ref{lemma 3.2} and the assumptions of the theorem. 
This shows that $P(\{R\neq M_{0}\})\rightarrow 0$. Following (\ref{buneal2.1}) with the choice of $q$, we can to conclude that $\hat{I}$ is consistent, i.e., $\displaystyle\lim_{n\rightarrow \infty} P(\widehat{I}=I_{0})=1$.\hspace{1cm} $\square$

\end{document}